\documentclass{amsart}
\usepackage{amsfonts}
\usepackage{amssymb}

\newtheorem{claim}{}[section]
\newtheorem{theorem}[claim]{Theorem}

\newtheorem{corollary}[claim]{Corollary}

\def\proclaim #1. #2\par{\medbreak
\noindent{\bf#1.\enspace}{\sl#2}\par\medbreak}

\def\proclaim #1. #2\par{\medbreak
\noindent{\bf#1.\enspace}{\sl#2}\par\medbreak}

\begin{document}
\large

\title[$q$-Wick products]{Feynman diagrams and Wick products \\associated with $q$-Fock
space}
\date{February 27, 2003}
\author{Edward G. Effros}
\address{Department of Mathematics\\
UCLA, Los Angeles, CA 90095-1555}
\email[Edward G. Effros]{ege@math.ucla.edu} 
\author{Mihai Popa}
\address{Department of Mathematics\\
UCLA, Los Angeles, CA 90095-1555}
\email[Mihai Popa]{mvpopa.ucla.edu}
\thanks{E. Effros was partially supported by the National Science
Foundation}

\begin{abstract}
It is shown that if one keeps track of crossings, Feynman diagrams
can be used to compute
$q$-Wick products and normal products in terms of each other.  
\end{abstract}

\maketitle

\setlength{\unitlength}{1cm}

\section{Introduction}

A recurrent theme in non-commutative analysis is that one may use graphs to
efficiently index the terms in complicated sums. One of the first to recognize this
principle was Cayley \cite{C1}, who introduced rooted trees
in order to label differentials (see \cite{Br} for additional examples). Currently,
the best-known example of graph-theoretic indexing may be found in perturbative
quantum field theory. In this context one uses Feynman diagrams to index
summands that arise when one evaluates the expectations of products of jointly
Gaussian random variables (see \cite{GJ}, \cite{Ja}). 

The random
variables of quantum field theory correspond to certain self-adjoint operators on
symmetric or antisymmetric Fock spaces (see \cite{GJ},\cite{Si}). In 1991, Bo\.{z}eko
and Speicher introduced a remarkable $q$-version of the Fock space, which for
$q=1,-1,$ and $0$ coincides with the symmetric (Boson), antisymmetric (Fermion), and
full (Voiculescu) Fock spaces (some of these ideas had been considered in
\cite{FB}; see \cite{BKS}). Bo\.{z}ejko and Speicher's
$q$-versions of stochastic processes and second quantization have attracted the
attention of a large number of researchers.

In \cite{BKS}, Bo\.{z}ejko, K\"{u}mmerer and Speicher introduced the $q$-analogs of
the Wick product. We show that some of the basic combinatorial
calculations involving Wick products of Gaussian random variables have natural
$q$-versions. In particular we use Feynman diagrams to express the Wick products in
terms of the normal (operator) products, and vice versa. We note that some special
cases of this material had been considered by Michael Anshelevich in \cite{An}.

We will explore $q$-forms of the Hopf algebraic theory of Kreimer and
Connes and Kreimer in a subsequent paper (see \cite{Kr}, \cite{CK}).

\section{$q$-Fock spaces and Feynman diagrams}

We begin by recalling the Bo\.{z}ejko and Speicher theory. Let $H_{0}$ be a
real Hilbert space and let $H$ be its complexification. We write 
\[
H^{\otimes n}=H\otimes \cdots \otimes H
\]
for the algebraic tensor product, and we define a hermitian form on the
algebraic sum 
\[
\mathcal{F}_{q}^{\rm alg}(H)=\Bbb{C}\Omega \oplus H\oplus H^{\otimes 2}\oplus
\cdots ,
\]
where $\Omega $ is taken to be a unit vector, by letting 
\[
\langle f_{1}\otimes \ldots \otimes f_{m},g_{1}\otimes \ldots \otimes
g_{n}\rangle _{q}=\delta _{mn}\sum_{\pi \in S_{n}}\langle f_{1},g_{\pi %
(1)}\rangle \ldots \langle f_{n},g_{\pi (n)}\rangle q^{\iota (\pi )},
\]
and where 
\[
\iota (\pi )=\#\left\{ (i,j):1\leq i<j\leq n,\pi (i)>\pi (j)\right\} 
\]
(we use $\#$ to indicate cardinality). We will generally delete the subscript $q$
in the hermitian form. The
\emph{$q$-Fock space}
$\mathcal{F}_{q}(H)
$ is the completion of this pre-Hilbert space. If $q=1$ or $q=-1,$ one must
first divide out by a null space, and one obtains the usual symmetric and
antisymmetric Fock spaces. If $q=0$ one obtains the usual full Fock space.
\emph{We restrict our attention to the case $-1<q<1.$} We let $H^{[0]}={\Bbb C}I$ and
\[H^{[n]}=\Bbb{C}\Omega \oplus H\cdots\oplus H^{\otimes n}.
\]
By an \emph{elementary tensor} we mean an element in $H^{\otimes n}$ of the from
$f_{1}\otimes\cdots \otimes f_{n}$.

For each $f\in H_{0}$, the creation and annihilation operators $a^{+}(f)$ and 
$a^{-}(f)$ are defined on $\mathcal{F}_{q}^{\rm alg}(H)$ by 
\begin{eqnarray}
a^{+}(f)\Omega  &=&f,  \nonumber \\
a^{+}(f)(f_{1}\otimes \cdots \otimes f_{n}) &=&f\otimes f_{1}\otimes \cdots
\otimes f_{n}  \label{creation}
\end{eqnarray}
and 
\begin{eqnarray}
a^{-}(f)\Omega  &=&0,  \nonumber \\
a^{-}(f)(f_{1}\otimes \cdots \otimes f_{n}) &=&\sum_{i=1}^{n}q^{i-1}\langle
f,f_{i}\rangle f_{1}\otimes \cdots \otimes \widehat{f_{i}}\otimes \ldots
\otimes f_{n}.  \label{annihilation}
\end{eqnarray}
$a^{+}(f)$ and $a^{-}(f)$ extend to bounded operators on $\mathcal{F}_{q}(H)$
satisfying
$a^{+}(f)=a^{-}(f)^{*}$. They satisfy the commutation relation
\[
a^{-}(f)a^{+}(g)-qa^{+}(g)a^{-}(f)=\langle f,g\rangle I.
\]

By analogy with the Boson case, we refer to the operators  
\[
\phi (f)=a^{+}(f)+a^{-}(f)
\]
as \emph{$q$-Gaussians}, and we define $\Gamma _{q}(H_{0})$ to be the von Neumann
algebra on
$\mathcal{F} _{q}(H)$ generated by these operators. We have that $\Omega $ is a
separating and cyclic vector for $\Gamma _{q}(H_{0})$. We let 
\[
\Bbb{E}:\Gamma _{q}(H_{0})\rightarrow \Bbb{C}:b\mapsto \langle b\Omega ,%
\Omega \rangle 
\]
be the corresponding state on $\Gamma _{q}(H_{0}).$ In particular if $\xi=\phi(f)$
and $\eta =\phi (g)$ for $f,g\in H_{0}$ then we have the \emph{covariance}
\[
\Bbb{E}(\xi \eta)=\langle a^{-}(f)a^{+}(g)\Omega,\Omega\rangle=\langle g,f
\rangle=\langle f,g \rangle
\]
(we recall that $H_{0}$ is a real Hilbert space).

We wish to compute the ``multivariable moments''
\[
\Bbb{E}(\xi_{1}\ldots\xi_{m})
\]
of $q$-Gaussians $\xi_{1},\ldots,\xi_{m}$. As in the classical case, these are
determined by polynomials of the ``covariances'' $\Bbb{E}(\xi_{i}\xi_{j})$, although
in this case the polynomials are non-commutative. 

We begin by letting
$a^{+1}(f)=a^{+}(f)$ and
$a^{-1}(f)=a^{-}(f).$ Our first task is to compute expressions of the form 
\[
m(\varepsilon )=\langle a^{\varepsilon (1)}(f_{1})a^{\varepsilon
(2)}(f_{2})\cdots a^{\varepsilon (2n)}(f_{2n})\Omega ,\Omega \rangle 
\]
for sequences $\varepsilon =(\varepsilon (1),\ldots ,\varepsilon (2n))\in
I(2n)=\left\{ 1,-1\right\} ^{2n}.$

Given $\varepsilon \in I(2n),$ we define 
\begin{eqnarray*}
\sigma _{2n} &=&\varepsilon (2n) \\
\sigma _{2n-1} &=&\varepsilon (2n-1)+\varepsilon (2n) \\
&&\ldots \\
\sigma _{2} &=&\varepsilon (2)+\cdots +\varepsilon (2n)\geq 0 \\
\sigma _{1} &=&\varepsilon (1)+\varepsilon (2)+\cdots +\varepsilon (2n).
\end{eqnarray*}
A simple induction shows that if $\sigma _{1}>0$, $\sigma _{2},\ldots
,\sigma _{k}\geq 0,$ then 
\begin{equation}\label{filter}
a^{\varepsilon (k)}(f_{k})\cdots a^{\varepsilon (2n)}(f_{2n})\Omega \in
H^{[\sigma _{k}]}, 
\end{equation}
and otherwise 
$a^{\varepsilon (k)}(f_{k})\cdots a^{\varepsilon (2n)}(f_{n})\Omega =0.$ We
say that $\varepsilon \in I(2n)$ is a \emph{Catalan sequence} if $\sigma
_{2n}>0,$ $\sigma _{2n-1}\geq 0,$ $\ldots ,\sigma _{1}=0,$ and we let $C(2n)$
be the set of such sequences. It is evident that $m(\varepsilon )=0$ unless $%
\varepsilon \in C(2n).$ We may associate a \emph{Catalan diagram} with each $%
\varepsilon \in C(2n)$. For example if $\varepsilon
=(-1,-1,1,-1,-1,1,1,1)\in C(8)$, then the corresponding Catalan diagram is
given by

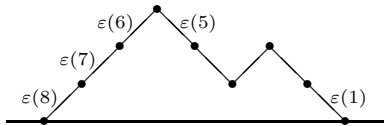
\[
\begin{picture}(6,2)
\put(0,.25){\circle*{.1}}
\put(0,.25){\line(1,1){.5}}
\put(.5,.75){\line(1,1){.5}}
\put(.5,.75){\circle*{.1}}
\put(1,1.25){\circle*{.1}}
\put(1,1.25){\line(1,1){.5}}
\put(1.5,1.75){\circle*{.1}}
\put(1.5,1.75){\line(1,-1){.5}}
\put(2,1.25){\circle*{.1}}
\put(2,1.25){\line(1,-1){.5}}
\put(2.5,.75){\circle*{.1}}
\put(2.5,.75){\line(1,1){.5}}
\put(3,1.25){\circle*{.1}}
\put(3,1.25){\line(1,-1){.5}}
\put(3.5,.75){\circle*{.1}}
\put(3.5,.75){\line(1,-1){.5}}
\put(4,.25){\circle*{.1}}
\put(-.5,.25){\line(1,0){5}}
\put(3.8,.5){\scriptsize$\varepsilon(1)$}
\put(.2,1){\scriptsize$\varepsilon(7)$}
\put(.7,1.5){\scriptsize$\varepsilon(6)$}
\put(1.8,1.5){\scriptsize$\varepsilon(5)$}
\put(-.3,.5){\scriptsize$\varepsilon(8)$}
\end{picture}
\]
where each $\varepsilon (k)$ is the slope of the corresponding line
segment. If $\varepsilon \in C(2n),$ we may evaluate $m(\varepsilon )$ in
terms of ``Feynman diagrams'' on $[2n]=\left\{ 1,\ldots ,2n\right\}$ (see \cite{GJ}
and \cite{Ja} for this terminology).

Given a finite linearly ordered set $S$, a \emph{Feynman diagram} $\gamma $
on $S$ is a partition of $S$ into one and two element sets. It will be
convenient to regard $\gamma \ $as a set of ordered pairs $\left\{
(i_{1},j_{1}),\ldots (i_{p},j_{p})\right\} $ with $i_{k}<j_{k}$ and $%
i_{k}\neq i_{l}$ and $j_{k}\neq j_{l}$ for $k\neq l$, and we refer to the unpaired
indices as {\em singletons}. With this notation we will assume that 
\begin{eqnarray}
i_{1}<\ldots <i_{n}  \label{order}
\end{eqnarray}
(the $j_{k}$ will generally be out of order). We write
$S^{+}$ (respectively, $S^{-}$) for the $j\in S$ with $(k,j)\in \gamma $ for some $k$
(respectively, the $i\in S$ such that $(i,k)\in \gamma $ for some $i$). We refer to
the elements of $S^{+}$ and $S^{-}$ as creators and annihilators, respectively. We may
specify a Feynman diagram with a simple graph of the form 
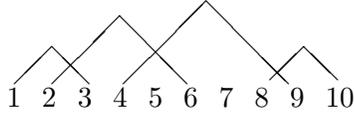
\[
\vspace{0.4cm}%
\begin{picture}(7,1.5)
\put(.1,.3){\line(1,1){.5}}
\put(.6,.8){\line(1,-1){.5}}
\put(.6,.3){\line(1,1){.9}}
\put(1.5,1.2){\line(1,-1){.9}}
\put(1.55,.3){\line(1,1){1.1}}
\put(2.65,1.4){\line(1,-1){1.1}}
\put(3.5,.35){\line(1,1){.45}}
\put(3.95,.8){\line(1,-1){.45}}
\put(0,0){1\hspace{.15cm} 2\hspace{.15cm} 3\hspace{.15cm} 4\hspace{.15cm}
5\hspace{.15cm} 6\hspace{.15cm} 7\hspace{.15cm} 8\hspace{.15cm} 9\hspace{.15cm}
10\hspace{.15cm}}
\end{picture}
\]

\noindent This corresponds to the Feynman diagram $\gamma
=\{(1,3),(2,6),(4,9),(8,10)\}$ on $[10]$. We let ${\mathsf F}(S)$ denote the set of
all Feynman diagrams on $S$, and we let ${\mathsf F}(n)={\mathsf F}([n])$.

We note
that more general partitions and their crossings are analyzed by using a
succession of semicircles to link  the elements of an equivalence class (see
\cite{Bi}).

We call the elements in $S$ the \emph{vertices} of the diagram. We say that
a pair $(k,l)\in \gamma $ is a \emph{left crossing }for $(i,j)$ if $k<i<l<j,$
and we define $c_{l}(i,j)$ to be the number of such left crossings. We refer
to $c(\gamma )=\sum_{(i,j)\in \gamma }c_{l}(i,j)$ as the \emph{crossing number}
of $\gamma $ (Biane calls this the ``restricted crossing number'' \cite{Bi}). The
total left crossings can be found by counting the intersections in the
corresponding graph. For example, in the above diagram, we have that
$c_{l}(1,3)=0$ and $c_{l}(2,6)=c_{l}(4,9)=c_{l}(8,10)=1$, and thus 
\[
c(\gamma )=c_{l}(1,3)+c_{l}(2,6)+c_{l}(4,9)+c_{l}(8,10)=3. 
\]
The general result is evident if one notes that an intersection in the graph
will occur between an ascending line for one pair $(k,l)\in \gamma $ and a
descending line for another $(i,j)\in \gamma $, which in turn will correspond
to the left crossing $(k,l)$ of $(i,j)$. We similarly have that $c(\gamma
)=\sum_{(i,j)\in \gamma }c_{r}(i,j)$, where $c_{r}(i,j)$ is the number of right
crossings $i<k<j<l$ where $(k,l)\in \gamma$. 

If $(i,j)\in \gamma ,$ we define
the \emph{gap} $g(i,j)$ to be the number of $k$ with $i<k<j$ and we let $%
a(i,j)=g(i,j)-c_{l}(i,j)$ (we will not need the right version of this). We define
$g(\gamma )=\sum_{(i,j)\in
\gamma }g(i,j)$ and $a(\gamma )=\sum_{(i,j)\in \gamma } a(i,j)=g(\gamma
)-c(\gamma ).$ Given
$(i,j)\in \gamma$ it is evident that $c_{l}(i,j)+c_{r}(i,j)\leq g(i,j)$,
and thus
\begin{eqnarray}
2c(\gamma)\leq g(\gamma)\label{inequality}
\end{eqnarray}

For some purposes, it is also useful to count \emph{degenerate crossings}.
These are the triples $i<k<j$, where $k$ is not paired and $(i,j)\in \gamma $%
. We let $d(\gamma)$ be the number of such triples in $\gamma$, and we
define the \emph{total crossing number} to be $tc(\gamma)=c(\gamma)+d(%
\gamma) $.

A Feynman diagram $\gamma $ on $S$ is \emph{complete} if the
there are no singletons (in which case $S$ must have an even number
of elements), and we let ${\mathsf F}_{c}(S)$ be the collection of all such diagrams.
Given
$\varepsilon
\in C(2n),$ we say that a complete Feynman diagram 
\[
\gamma =\left\{ (i_{1},j_{1}),\ldots ,(i_{n},j_{n})\right\} 
\]
on $[2n]$ is \emph{compatible} with $\varepsilon $ if $\varepsilon (i_{k})=-1$
and $\varepsilon (j_{k})=1.$ We let $\mathsf {F_{c}}(\varepsilon )$ be the set of all
such Feynman diagrams on $S$. It is easy to see that $\mathsf {F_{c}}(\varepsilon )$
is non-empty. For example, we may fix a vertex of maximum height, and then pair the
decending and ascending edges adjacent to that vertex. In the above diagram we begin
by placing the pair $(5,6)$ in $\gamma $. Eliminating the two line segments
and rejoinging the graph, we can continue by induction to appropriately pair
all the ascending edges with descending edges.  Conversely, each complete Feynman
diagram $\gamma $ on $[2n]$ is compatible with the unique sequence $\varepsilon \in
C(2n)$ defined by letting $\varepsilon (k)=1$ if $k\in [2n]^{+}$ and $\varepsilon
(k)=-1$ if $ k\in [2n]^{-}$.

Given $q$-Gaussian random variables $\xi _{i}=\phi (f_{i})$ ($i=1,\ldots ,n)$
and a Feynman diagram $\gamma $ on $[n]$, we let 
\begin{eqnarray*}
v(\gamma )
&=&\Bbb{E}(\xi_{i_{1}}\xi_{j_{1}})\cdots \Bbb{E}(
\xi_{i_{p}}\xi_{j_{p}}) \xi _{h_{1}}\ldots \xi _{h_{r}}
\end{eqnarray*}
where $\gamma =\left\{ (i_{k},j_{k}):k=1,\ldots ,p\right\} ,$and $%
h_{1}<h_{2}<\ldots <h_{r}$ are the $\gamma$ singletons. 

\begin{theorem}
\label{t1} For any $\varepsilon \in C(2n),$ we have that 
\[
m(\varepsilon )=\sum_{\gamma \in \mathsf{ F_{c}}(\varepsilon )}v(\gamma )q^{c(\gamma
)}.
\]
\end{theorem}

\proof Given $\gamma \in {\mathsf F}_{c}(\varepsilon ),$ we will as before assume
that $\gamma=\{(i_{n},j_{n})\}$ where $i_{1}<\ldots<i_{n}$. We define a sequence of
elementary tensors 
\[\Omega _{2n}^{\gamma },\ldots ,\Omega _{0}^{\gamma }
\]
 as follows. Let us suppose that $i_{n}<i_{n}+1<\ldots <j_{n}\leq 2n.$ We define 
\begin{eqnarray*}
\Omega _{2n}^{\gamma } &=&a^{+}(f_{2n})\Omega =f_{2n} \\
\Omega _{2n-1}^{\gamma } &=&a^{+}(f_{2n-1})\Omega _{2n}^{\gamma
}=f_{2n-1}\otimes f_{2n} \\
&&\cdots  \\
\Omega _{i_{n}+1}^{\gamma } &=&a^{+}(f_{i_{n}+1})\Omega %
_{2n-k+1}=f_{i_{n}+1}\otimes f_{i_{n}+2}\otimes \cdots \otimes
f_{j_{n}-1}\otimes f_{j_{n}}\otimes\cdots
\end{eqnarray*}
From (\ref{annihilation}), 
\begin{eqnarray*}
a^{-}(f_{i_{n}})(\Omega _{i_{n}+1}^{\gamma }) &=&\langle
f_{i_{n}},f_{i_{n+1}}\rangle \widehat{f_{i_{n}+1}}\otimes f_{i_{n}+2}\otimes
\cdots +\cdots  \\
&&+q^{j_{n}-i_{n}-1}\langle f_{i_{n}},f_{j_{n}}\rangle f_{i_{n}+1}\otimes
\cdots \otimes \widehat{f_{j_{n}}}\otimes \cdots +\cdots 
\end{eqnarray*}
We define $ \Omega _{i_{n}}^{\gamma }$ to be a particular elementary tensor summand in
this expression: 
\[
\Omega _{i_{n}}^{\gamma }=q^{j_{n}-i_{n}-1}\langle
f_{i_{n}},f_{j_{n}}\rangle f_{i_{n}+1}\otimes \cdots \otimes \widehat{%
f_{j_{n}}}\otimes \cdots .
\]
If $i_{n}<k<j_{n},$ then we must have that $k\in [2n]^{+},$ i.e., there is
an $h$ with $(h,k)\in \gamma .$ Since $h<i_{n},$ it follows that $(h,k)$ is
a left crossing for $(i_{n},j_{n})$ and we see that there are exactly $%
c(i_{n},j_{n})=j_{n}-i_{n}-1$ terms between $i_{n}$ and $j_{n}.$ We conclude
that 
\[
\Omega _{i_{n}}^{\gamma }=q^{c(i_{n},j_{n})}\langle
f_{i_{n}},f_{j_{n}}\rangle f_{i_{n}+1}\otimes \cdots \otimes \widehat{%
f_{j_{n}}}\otimes \cdots .
\]

Let us suppose that we have defined elementary tensors $\Omega _{2n}^{\gamma
},\Omega _{2n-1}^{\gamma },\ldots ,\Omega _{k}^{\gamma }$ in such a manner
that $i_{p}<k\leq i_{p+1},$ and none of the factors $f_{j_{p}+1},\ldots
,f_{j_{n}}$ occurs in $\Omega _{k}^{\gamma }$. If $k>i_{p}+1$ , we let 
\[
\Omega _{k-1}^{\gamma }=a^{+}(f_{k})\Omega _{k}^{\gamma }=f_{k}\otimes\Omega
_{k}^{\gamma }.
\]
 On the other
hand if $k=i_{p}+1,$ let us suppose that 
\[
\Omega _{k}^{\gamma }=\alpha (f_{k_{1}}\otimes \cdots \otimes
f_{k_{r}}\otimes f_{j_{p}}\otimes f_{\ell _{1}}\otimes \ldots ). 
\]
for some scalar $\alpha .$ We define 
\[
\Omega _{i_{p}}^{\gamma }=\Omega _{k-1}^{\gamma }=q^{r}\alpha \langle
f_{i_{p}},f_{j_{p}}\rangle f_{k_{1}}\otimes \ldots \otimes f_{k_{r}}\otimes 
\widehat{f_{j_{p}}}\otimes f_{\ell _{1}}\otimes \ldots , 
\]
which is one of the elementary tensor summands of $a^{-}(f_{i_{p}})\Omega
_{i_{p}+1}^{\gamma }.$ Each $k_{t}$ lies in $[2n]^{+}\setminus
\{j_{p+1},\ldots,j_{n}\}$ and $k_{t}<j_{n}$, hence if
$k_{t}=j_{h},$ then $h<p$ and therefore $i_{h}<i_{p}.$ It follows that
$ (i_{h},k_{t})$ is a left crossing for $(i_{p},j_{p}).$ Since this holds for each
$t$ we see that
$r=c(i_{p},j_{p}).$
Continuing in this manner $\Omega _{0}^{\gamma }=v(\gamma
)q^{c(\gamma )}\Omega .$ 

Since it is evident that every non-zero elementary tensor summand in  
\[
a^{\varepsilon (1)}(f_{1})\cdots a^{\varepsilon (2n)}(f_{2n})\Omega
\]
corresponds to a unique Feynman diagram in $\mathsf {F_{c}}(\varepsilon)$,
\[
m(\varepsilon)=\sum_{\gamma \in \mathsf {F_{c}}(\varepsilon)} \langle\Omega
_{0}^{\gamma },\Omega\rangle=\sum_{\gamma \in \mathsf
{F_{c}}(\varepsilon)}v(\gamma)q^{c(\gamma )},
\]
and we are done.
 \endproof

The following is due to Anshelevich \cite{An}.

\begin{corollary}[$q$-Wick theorem]\label{qwick}
For any $q$-Gaussian random variables $\xi _{i},$ we have 
\[
\Bbb{E}(\xi _{1}\xi _{2}\ldots \xi _{2n})=\sum_{\gamma \in \mathsf {F_{c}}(2n)}
v(\gamma )q^{c(\gamma )}
\]
and on the other hand, 
\[
\Bbb{E}(\xi _{1}\xi _{2}\ldots \xi _{2n+1})=0.
\]
\end{corollary}

\proof We have that 
\begin{eqnarray*}
\Bbb{E}(\xi _{1}\ldots \xi _{2n}) &=&\langle (a^{-}(f_{1})+a^{+}(f_{1}))\cdots
(a^{-}(f_{2n})+a^{+}(f_{2n}))\Omega ,\Omega \rangle \\
&=&\sum_{\varepsilon \in I(2n)}m(\varepsilon )=\sum_{\varepsilon \in
C(n)}m(\varepsilon ),
\end{eqnarray*}
and thus the formula follows from the theorem. The result for odd moments is
immediate from (\ref{filter}). \endproof

\section{$q-$Wick products}

Following \cite{BKS}, we define the $q$-\emph{Wick product} for $f_{1},\ldots
,f_{n}\in H$ by 
\begin{equation}
W(f_{1},\ldots ,f_{n})=\sum a^{+}(f_{i_{1}})\cdots
a^{+}(f_{i_{k}})a^{-}(f_{j_{1}})a^{-}(f_{j_{2}})\ldots a^{-}(f_{j_{l}})q^{\iota (I,J)}
\label{Wick1}
\end{equation}
where the sum is taken over all families $I=\left\{ i_{1},\ldots
,i_{k}\right\} $ and $J=\left\{ j_{1},\ldots ,j_{l}\right\} $ where $%
i_{1}<\ldots <i_{k},$ $j_{1}<\ldots <j_{l},$ and $I\sqcup J=[n],$ and we let 
$\iota (I,J)=\#(\left\{ (p,q):i_{p}>j_{q}\right\} .$ This operator is
characterized by the recursion 
\begin{eqnarray}
W(f,f_{1},\ldots f_{n})&=&\phi (f)W(f_{1},\ldots ,f_{n})\nonumber\\
&&\hspace{.4cm}-\sum_{i=1}^{n} q^{i-1}\langle
f,f_{i}\rangle W(f_{1},\ldots \hat{f}_{i},\ldots ,,f_{n})  \label{Wickind}
\end{eqnarray}
(see \cite{BKS}, proof of Prop. 2.7). It follows from a simple induction on (\ref
{Wickind}), that $W(f_{1},\ldots f_{n})$ is a non-com\-muting polynomial of
the operators $\xi _{i}=\phi (f_{i})$, $(1\leq i\leq n)$, and we will use
the usual Wick product notation 
\[
:\!\xi _{1}\ldots \xi _{n}\!:=W(f_{1},\ldots ,f_{n}). 
\]
Since $\Omega $ is separating for $\Gamma _{q}(H_{0}),$ $u=:\!\xi _{1}\ldots
\xi _{n}\!:$ is the unique operator in $\Gamma _{q}(H_{0})$ satisfying 
\[
u\Omega =f_{1}\otimes \ldots \otimes f_{n}. 
\]

The following provides an explicit (non-recursive) expression for the $q$%
-Wick product analogous to the ``classical formula '' for $q=1$ (see \cite{Ja},
Th. 3.4). We let $\#(\gamma )$ denote the number of pairs in $\gamma .$

\begin{theorem}
For any $q$-Gaussian random variables $\xi _{i}=\phi (f_{i})$ ($i=1,\ldots ,n)$
\begin{equation}
:\!\xi _{1}\cdots \xi _{n}\!:=\sum_{\gamma \in {\mathsf F}(n)} (-1)^{\#(\gamma
)}q^{a(\gamma )}v(\gamma ).\label{Wick2}
\end{equation}
\end{theorem}

\proof ${\mathsf F}(1)$ contains only the empty Feynman diagram
$\gamma_{0}$, and  
\[
:\!\xi_{1}\!:=W(f_{1})=a^{+}(f_{1})+a^{-}(f_{1})=\xi_{1}=v(\gamma_{0}).
\]
On the other hand, ${\mathsf F}(2)=\{\gamma_{0},\gamma_{1}\}$ where $\gamma_{0}$ is
empty and $\gamma_{1}=\{(1,2)\}$. We have
\begin{eqnarray*}
:\!\xi_{1}\xi_{2}\!:&=& W(f_{1},f_{2})\\
&=&a^{+}(f_{1})a^{+}(f_{2})+a^{+}(f_{1})a^{-}(f_{2})\\&&\hspace{.2cm}+qa^{+}(f_{2})a^{-}(f_{1})+a^{-}(f_{1})a^{-}(f_{2})\\
&=&(a^{+}(f_{1})+a^{-}(f_{1}))(a^{+}(f_{2})+a^{-}(f_{2}))-\Bbb{E}(\xi_{1}\xi_{2})I\\
&=&\xi_{1}\xi_{2}-\Bbb{E}(\xi_{1}\xi_{2})I\\
&=&v(\gamma_{0})-v(\gamma_{1}).
\end{eqnarray*}

Let us suppose that we have proved (\ref{Wick2}) for $n-1$ and $n$. If $\xi _{0}=\phi
(f_{0}),$ then applying the formula for $n$ and $n-1$ and the recurrence relation, 
\begin{eqnarray*}
W(f_{0},f_{1},\ldots f_{n})&=&\xi _{0}\sum_{\gamma \in {\mathsf F}(n)}(-1)^{\#(\gamma
)}q^{a(\gamma )}v(\gamma )\\
&&-\sum_{l=1}^{n}q^{l-1}\langle f_{0},f_{l}\rangle 
\sum_{\delta \in {\mathsf F}([n] \setminus\{l\})} (-1)^{\#(\delta
)}q^{a(\delta )}v(\delta)
\end{eqnarray*}

Each Feynman diagram $\gamma =\left\{
(i_{k},j_{k})\right\} \in {\mathsf F}(n)$ trivially determines a Feynman diagram $%
\gamma ^{\prime }=\left\{ (i_{k},j_{k})\right\} \in {\mathsf F}(\{ 0\}\cup[n]\}
),$ for which $\xi_{0}v(\gamma)=v(\gamma^{\prime})$. Since $\#(\gamma ^{\prime
})=\#(\gamma )$ and
$a(\gamma ^{\prime })=a(\gamma ),$  
\[
\xi_{0} (-1)^{\#(\gamma )}q^{a(\gamma )}v(\gamma )=
(-1)^{\#(\gamma ^{\prime })}q^{a(\gamma ^{\prime })}v(\gamma ^{\prime
}).
\]

Each Feynman diagram $\delta \in {\mathsf F}([n] \setminus\{l\})$
determines a Feynman diagram $\delta ^{\prime }=\delta\cup \left\{ (0,l)\right\}
\in {\mathsf F}(\{ 0\}\cup[n] )$ for which 
\[
\langle f_{0},f_{l}\rangle v(\delta)=v(\delta^{\prime})
\]
It is evident that $\#(\delta
^{\prime })=\#(\delta )+1.$ We have that $c_{\delta ^{\prime }}(0,l)=0$ and
thus $a_{\delta ^{\prime }}(0,l)=l-1.$ If $(i,j)\in \delta ,$ then we may
consider three cases. If $l<i,$ or $j<l,$ then it is evident that $a_{\delta
^{\prime }}(i,j)=a_{\delta }(i,j)$. If $i<l<j,$ then $g_{\delta
^{\prime }}(i,j)=g_{\delta }(i,j)+1.$ On the other hand $0<i<l<j$
introduces another left crossing when we regard $(i,j)$ as an element of $%
\delta ^{\prime },$ and thus $c_{\delta ^{\prime }}(i,j)=c_{\delta }(i,j)+1.$
It follows that $a(\delta ^{\prime })=l-1+a(\delta )$ and 
\[
-q^{l-1}\langle f_{0},f_{l}\rangle (-1)^{\#(\delta
)}q^{a(\delta )}v(\delta)=(-1)^{\#(\delta ^{\prime })}q^{a(\delta
^{\prime })}v(\delta ^{\prime }). 
\]

A Feynman diagram $\theta $ in ${\mathsf F}(\{ 0\}\cup[n] )$ has the form
$\gamma^{\prime}$ if and only if ${0}$ is a singleton in $\theta$, and the form
$\delta^{\prime}$ if  $(0,l)\in \theta$. Thus
\begin{equation}
:\!\xi _{0}\xi _{1}\cdots \xi _{n}\!:=\sum_{\theta \in {\mathsf F}(\{ 0\}\cup[n])} 
(-1)^{\#(\theta )}q^{a(\theta )}v(\theta ),
\end{equation}
and we are done.

\endproof

Conversely, we may express products $\xi _{1}\ldots \xi _{n}$ in terms of $q$%
-Wick products. For this purpose we need to generalize the $q$-Wick theorem to
products of $q$-Wick products. Given  
$q$-Gaussian random variables $\{\xi _{p,k}\}$ with $1\leq p
\leq t$ and $1\leq k \leq n_{p}$, we may regard the index set $S=\{(p,k)\}$ as
partitioned by the first integer, and we refer to each partition as a {\em block}. We
let  
$S_{\rm lex}$ denote $S$ with the lexicographic ordering 
\[
(1,1)<(1,2)<\cdots (2,1)<\cdots (t,n_{t})
\]

\begin{theorem}\label{expectprod}
Suppose that we are given Gaussian random variables $\{\xi _{p,k}\}$ with $1\leq p
\leq t$ and $1\leq k \leq n_{p}$, 
Then if  $Y_{p}=:\!\xi_{p,1}\cdots \xi _{p,n_{p}}\!:$, we have 
\[
\Bbb{E}(Y_{1}\cdots Y_{t})=\sum v(\gamma )q^{c(\gamma )}
\]
where the sum is taken over all complete Feynman diagrams $\gamma $ on $S_{\rm lex}$
which do not link vertices within blocks.
\end{theorem}
\proof A typical summand of $Y_{1}\cdots Y_{t}$ has the form 
\[
u=\prod_{p=1}^{t}a^{+}(f_{i_{1}}^{(p)})\cdots
a^{+}(f_{i_{r(p)}}^{(p)})a^{-}(f_{j_{1}}^{(p)})a^{-}(f_{j_{2}}^{(p)})\ldots
a^{-}(f_{j_{s(p)}}^{(p)})q^{\iota (I_{p},J_{p})}. 
\]
where for each $p,$ $j_{1}<j_{2}<\ldots $ and $i_{1}<i_{2}<\ldots $.

We may use Theorem
\ref{t1} to compute $\langle u\Omega,\Omega \rangle$  provided we reorder the index
set. We let $S_{u}$ denote $S$ with the total ordering $(p,k)<(p^{\prime
},k^{\prime })$ if $p<p^{\prime },$ and 
\begin{equation}
(p,i_{1})<(p,i_{2})<\ldots <(p,j_{1})<(p,j_{2})<\ldots \label{modorder}
\end{equation}
Thus if we let $\varepsilon (p,i_{g})=1$ and $\varepsilon (p,j_{h})=-1,$ 
\[
\langle u\Omega,\Omega \rangle=m(\varepsilon)=\sum_{\gamma \in
\mathsf{F_{c}}(\varepsilon )}v(\gamma )q^{c(\gamma )}q^{\iota (I,J)} 
\]
where $\iota (I,J)=\sum\iota(I_{k},J_{k})$. It should be noted that if $\gamma \in
\mathsf{F_{c}}(\varepsilon )$, then $%
\gamma $ will not link elements of a block, since in the definition of
$v(\gamma)$, a creator is always paired with an annihilator on its left. 

We may use a sequence of $\iota (I,J)$ transpositions of the index set
$S$ to transform $S_{u}$ into $S_{\rm lex}$. Retaining the same ordered pairs, each
Feynman diagram $\gamma$ on $S$ with a given total ordering is mapped  to a
Feynman diagram $\gamma^{\prime}$ on the reordered set. It is evident that
$v(\gamma)=v(\gamma^{\prime})$, but in general the number of crossings will change.

If $S_{u}\neq S_{\rm lex}$, we must have that $i_{r(p)}>j_{1}$ for some $p$. Our first
step will be the transition from 
\[
(p,i_{1})<\ldots <(p,i_{r(p)})<(p,j_{1})<\ldots <(p,j_{s(p)}) 
\]
to 
\[
(p,i_{1})<\ldots <(p,j_{1})<(p,i_{r(p)})<\ldots <(p,j_{s(p)}), 
\]
Continuing in this manner, a series of $\iota (I,J)$ transpositions
will give us a chain of Feynman diagrams $\gamma _{0}\rightarrow \gamma
_{2}\rightarrow \ldots \rightarrow \gamma _{\iota (I,J)}$ on the permuted sets
which will return us to the original ordering.

At each stage we will have an adjacent ``disordered'' pair $(p,j_{l}),(p,i_{k})$ with
$ j_{l}>i_{k}$ and we perform the transposition 
\[
(p,j_{l})\leftrightarrow (p,i_{k}). 
\]
Let us consider the crossings in the corresponding Feynman diagrams $\gamma $
and $\gamma ^{\prime }.$ Suppose that $a<b<c<d$ is a crossing in either $%
\gamma $ or $\bar{\gamma}.$ If none of these four vertices coincides with $%
i\,_{k}$ or $j_{l},$ the crossing will be left invariant under the
transpositions $\gamma \leftrightarrow \gamma ^{\prime }.$ We have three
other cases (remember that $c,d$ and $(p,j_{k})$ are creators, and the terms 
$(p,i_{k}),$ $(p,j_{l})$ are adjacent) 
\begin{eqnarray*}
a &<&b<(p,i_{k})<c=(p,j_{l})<d \\
a &<&b=(p,i_{k})<(p,j_{l})<c<d \\
a &<&b=(p,i_{k})<c=(p,j_{l})<d
\end{eqnarray*}
In the first and second cases, the crossing $a<b<c<d$ will remain unaffected
under these transpositions. From our construction, the third case will occur
precisely when the non-crossing sequence $a<(p,j_{l})<(p,i_{k})<d$ occurs in $%
\gamma $ with $j_{l}>i_{k}$, and we obtain $\gamma ^{\prime }$ by the
transposition It follows that $c(\gamma ^{\prime })=c(\gamma )+1,$ and thus $%
q^{c(\gamma )}q^{h}=q^{c(\gamma ^{\prime })}q^{h-1}.$ 

Starting with a Feynman
diagram $\gamma=\gamma_{0}$ on $S_{u}$, we obtain a non-crossing Feynman diagram
$\gamma ^{\prime }=\gamma_{\iota(I,J)}$ on
$S_{lex}$ with $v(\gamma)q^{c(\gamma )}q^{\iota(I,J)}=v(\gamma^{\prime})q^{c(\gamma
^{\prime })}$. It is easy to see that all complete diagrams on $S_{lex}$ which do
not link elements within blocks arise in this fashion, hence taking the sum of
terms $u$, we obtain the desired result.
\endproof

\begin{theorem}
Let $Y_{i}=:\!\xi _{i1}\cdots \xi _{in_{i}}\!:$ Then 
\begin{eqnarray}
Y_{1}\cdots Y_{t}=\sum :\!v(\gamma )\!:q^{tc(\gamma )}\label{prodwick}
\end{eqnarray}
where the sum is taken over all Feynman diagrams $\gamma $ on $S_{lex}$
which do not link vertices within blocks.
\end{theorem}
\proof Let $A=Y_{1}\cdots Y_{t}$ and $B=\sum :\!v(\gamma )\!:q^{tc(\gamma )}.$
As in the proof of \cite{Ja}, Th. 3.15, it suffices to show that for any Wick product
$W=:\!\eta _{1}\cdots
\eta _{u}\!:,$ $\Bbb{E}(AW)=$ $\Bbb{E}(BW).$ From Corollary \ref{qwick},
\[
\Bbb{E}(AW)=\sum v(\delta )q^{c(\delta )}, 
\]
where we sum over all complete Feynman diagrams $\delta $ on
the partitioned ordered set 
\[
T=\{(1,1),\ldots ,(1,n_{1});(2,1),\ldots ,(t,n_{t});1,\ldots ,u\} 
\]
(where the semicolons separate blocks) which do not link vertices in the same block.
Each such diagram determines a (generally incomplete) diagram $\gamma =\delta \cap
(S\times S)$ on $S_{lex}$ with the same property. Since $\delta $ is complete, the
singletons
$(p_{1},k_{1})<\ldots< (p_{u}, k_{u})$ in
$\gamma$ are linked with elements of $[u]=\{1,\ldots,u\}$. It follows that 
\begin{eqnarray*}
v(\delta )&=&\left (\prod_{(i,j)\in \gamma }\Bbb{E}( \xi _{i}\xi _{j})\right ) \Bbb{E}(
\xi _{(p_{1},k_{1})}\eta _{g_{1}}) \cdots \Bbb{E}( \xi
_{(p_{u},k_{u})}\eta _{g_{u}})\\
&=& \left (\prod_{(i,j)\in \gamma }\Bbb{E}( \xi _{i}\xi _{j})\right )v(\theta)
\end{eqnarray*}
where $\theta$ is a complete Feynman diagram on 
\[
T_{0}=\left\{ (p_{1},k_{1}),\ldots ,(p_{u},k_{u});1,\ldots ,u\right\} 
\]
which does not link vertices in the two blocks. 

On the other hand, given a Feynman diagram $\gamma $ on $S$ which
does not link elements in any of the blocks,  
\[
v(\gamma )=\left (\prod_{(i,j)\in \gamma }\Bbb{E}( \xi _{i}\xi _{j})\right ) \xi
_{(p_{1},k_{1})}\cdots \xi _{(p_{u},k_{u})} 
\]
and thus 
\[
:\!v(\gamma )\!:=\left (\prod_{(i,j)\in \gamma }\Bbb{E}( \xi _{i}\xi _{j})\right )
:\!\xi _{(p_{1},k_{1})}\cdots \xi _{(p_{u},k_{u})}\!:. 
\]
From Theorem \ref{expectprod}, 
\begin{eqnarray*}
\Bbb{E}(:\!v(\gamma )\!:W) &=&\left (\prod_{(i,j)\in \gamma }\Bbb{E}( \xi _{i}\xi
_{j})\right ) \Bbb{E}(:\!\xi _{(p_{1},k_{1})}\cdots \xi _{(p_{u},k_{u})}\!:W) \\
&=&\left (\prod_{(i,j)\in \gamma }\Bbb{E}( \xi _{i}\xi _{j})\right)
\sum_{\theta \in \mathsf{G}} v(\theta )q^{c(\theta )}
\end{eqnarray*}
where $\mathsf{G}$ is the set of all complete Feynman diagrams $\theta$ on $T_{0}$
which do not link vertices of the two blocks.
Given such a $\theta ,$ $\delta =\gamma \cup
\theta $ is a generic complete Feynman diagram on $T$ extending $\gamma 
$ which does not link internal vertices. It is evident that 
\[
c(\delta )=c(\gamma )+c(\theta)+d(\gamma ) 
\]
since forming the term $\prod_{(i,j)\in \gamma }\Bbb{E}( \xi
_{i}\xi _{j}) $ will ``hide'' precisely $d(\gamma )$ intersections
arriving from the pairs in $\theta.$ It follows that 
\[
\Bbb{E}(:\!v(\gamma )\!:W)q^{ct(\gamma )}=\prod_{(i,j)\in \gamma }\Bbb{E}( \xi
_{i}\xi _{j}) \sum_{\theta \in \mathsf{G}} v(\theta
)q^{c(\theta )+ct(\gamma )}=\sum v(\delta )q^{c(\delta )}, 
\]
where we sum over non-linking extensions $\delta$ of $\gamma$, and thus
\[
\Bbb{E}(BW)=\sum :\!v(\gamma )\!:q^{tc(\gamma )}=\sum v(\delta )q^{c(\delta
)}=\Bbb{E}(AW).
\]

\endproof

\begin{corollary} {\rm (See \cite{An} Remark 6.15 for a one variable version of
this result)}. Given
$\xi_{j}=\phi(f_{j})$ as above, we have
\[
\xi_{1}\cdots \xi_{n}=\sum_{\gamma \in \mathsf{F}(n)} :\!v(\gamma )\!:q^{tc(\gamma )}.
\]
\end{corollary}

\section{The free case ($q=0$) and non-crossing diagrams}

We say that a Feynman diagram $\gamma$ on a finite totally ordered set $S$ is \emph{
non-crossing} if $c(\gamma)=0$, \emph{\ strongly non-crossing} if $tcr(\gamma)=0$,
and that it is \emph{gap free} if
$g(\gamma)=0$. We let
$\mathsf{NC}(S)$, $\mathsf{SNC}(S)$, and $\mathsf{GF}(S)$ denote the
corresponding diagrams on
$S$. It is evident that 
\[
\mathsf{NC}(S)\supseteq \mathsf{SNC}(S)\supseteq \mathsf{GF}(S).
\] 
If $S=[n]$, we will simply write
$\mathsf{NC}(n)$, etc. 

If $q=0,$ we have the commutation relation $a^{-}(f)a^{+}(g)=<f,g>I$. The convention
of \cite{BS1} is that $q^{c}=0$ for
$c\neq 0,$ and $q^{0}=1.$ Thus we may drop terms in the $0$-Wick theorem for which $q$
is raised to a positive power:
\[
\Bbb{E}(\xi _{1}\xi _{2}\cdots \xi _{2n})=\sum_{\gamma \in
\mathsf{NC}(2n)}v(\gamma )
\]
 
Turning to Wick products, we delete terms with inversions in the 
definition of the free Wick product:
\[
W(f_{1},\ldots ,f_{n})=\sum_{k=0}^{n} a^{+}(f_{1})\cdots a^{+}(f_{k})
a^{-}(f_{k+1})\cdots a^{-}(f_{n}).
\]
For the alternative formula (\ref{Wick2}) we need only consider
terms with $a(\gamma )=g(\gamma )-c(\gamma )=0$.  It follows from (\ref{inequality})
that $g(\gamma)=0$, and thus $\gamma $ is gap free.  We conclude that 
\begin{equation}
:\!\xi _{1}\cdots \xi _{n}\!:=\sum_{\gamma \in \mathsf{GF}(n)} (-1)^{\#(\gamma
)}v(\gamma ).
\end{equation}
We have, for example, that 
\begin{equation}
:\!\xi _{1}\xi _{2}\xi_{3}\!:=\xi _{1}\xi _{2}\xi_{3}-{\Bbb E}(\xi _{1}\xi _{2})\xi
_{3} - {\Bbb E}(\xi _{2}\xi _{3}) \xi_1.
\end{equation}
Turning to (\ref{prodwick}), we have
\begin{eqnarray}
Y_{1}\cdots Y_{t}=\sum :\!v(\gamma )\!:
\end{eqnarray}
where the sum is over strongly non-crossing diagrams which do not link elements of a
block. Thus in particular,
\begin{equation}
\xi _{1}\cdots \xi _{n}=\sum_{\gamma \in \mathsf{SNC}(S)} :\!v(\gamma )\!:.
\end{equation}

\end{document}